\documentclass{statsoc}
\usepackage[dvips]{graphicx}
\usepackage{amsbsy, amssymb, amsmath}

\usepackage[usenames]{color}
\usepackage{colortbl}

\definecolor{light}{gray}{0.9}
\definecolor{dark}{gray}{0.8}
\definecolor{darkest}{gray}{0.7}

\newcommand{\be}{\begin{equation}}
\newcommand{\ee}{\end{equation}}
\newcommand{\bea}{\begin{eqnarray}}
\newcommand{\eea}{\end{eqnarray}}
\newcommand{\bbe}{\begin{eqnarray*}}
\newcommand{\ebe}{\end{eqnarray*}}
\newcommand{\bean}{\begin{eqnarray*}}
\newcommand{\eean}{\end{eqnarray*}}

\newcommand{\bm}{\boldmath}

\newtheorem{theorem}{Theorem}
\newtheorem{LEM}{Lemma}

\newtheorem{definition}{Definition}
\newtheorem{example}{Example}

\newcommand{\argmin}{\mbox{arg}\min}          

\renewcommand{\d}{\mbox{d}}

\newcommand{\RAIC}{ QAIC } 
\newcommand{\RBIC}{ QBIC } 
\newcommand{\RADE}{ MRA }

\newcommand{\0}{-}
\renewcommand{\t}[1]{{$ \mathbf{#1^*}$}}
\newcommand{\simtitle}[1]
{ \hline
\rowcolor[gray]{0.85}&&&&&&&&&&\\[-2.5ex]
\multicolumn{11}{>{\columncolor[gray]{0.85}}l }{{ #1}}\\
\hline}
    
\newcommand{\simtitleg}[1]
{ \hline
\rowcolor[gray]{0.85}&&&&&&&&&&&&&&&&\\[-.9ex]
\multicolumn{17}{>{\columncolor[gray]{0.85}}l }{{ #1}}\\[1ex]
\hline}

\newtheorem{PRO}{Proposition}[section]

\title[Quadratic Risk]{
Model selection in High-Dimensions: A Quadratic-risk based approach \\
}

\author[Ray and Lindsay]{ Surajit Ray }
\address{ Boston University, Boston,
  USA.}
\author[Ray and Lindsay]{Bruce G. Lindsay}
\address{ Pennsylvania State University, 
  University Park,
  USA.}
\email{sray@bios.unc.edu}

\begin{document}

\maketitle

\begin{abstract}
In this article we propose a general class of risk measures which can be
used for data based evaluation of parametric models.
The loss  function is defined as  generalized quadratic distance between the true
density and the proposed  model. These distances are
characterized by a simple quadratic form structure that is adaptable through
the choice of a nonnegative definite kernel and a bandwidth parameter. Using  asymptotic results for
the quadratic distances we build a quick-to-compute
approximation for the risk function. Its derivation is analogous to
the Akaike Information Criterion (AIC), but unlike AIC, the quadratic risk is a global comparison tool.
The method does not require resampling, a great advantage when point estimators are expensive to compute. The method is illustrated using the problem of selecting the number of components in a mixture model, where it is shown that, by using an appropriate kernel, the method is computationally straightforward in arbitrarily high data dimensions. In this same context it is shown that the method has some clear advantages over AIC and BIC. 


\end{abstract}
\medskip
\small
\textit{Address for correspondence:} Surajit Ray, Department of
Mathematics and Statistics,
 111 Cummington Street, Boston, MA 02215, USA\\
Email: sray@math.bu.edu

\normalsize

\newpage

\section{Introduction}

In this article we consider data-based evaluation of statistical 
models, where by a statistical model we mean a parametric family of
distributions. We will denote a  parametric model by $\mathcal{M}$, and an element of the model    by  $M_\theta$, for $\theta$ in some
parameter space $\Theta$, so that  $\mathcal{M}=\{M_\theta ;\theta \in \Theta\}$. A collection of
a set of probability models $\mathcal{M}_k$, $k=1,2,\ldots$ will constitute a
model class, say $\mathbb{M}$.
We will focus on the problem of selecting one
statistical model from  the competing models in the class $\mathbb{M}=\{\mathcal{M}_k\}$, 
where $k$ is an index for the parametric model under consideration. 
In the example we consider in Section~\ref{selection}, 
$\mathbb M$ is the class of all finite mixtures of multivariate
  Gaussians, and for each $k$, $\mathcal{M}_k$ is the set of all $k$-component
  multivariate Gaussian
  mixtures. A particular model element $M_\theta \in \mathcal{M}_k$
  is a
  $k$-component mixture with specific values for the parameters $\theta$.
Let us also use $F_\tau$ to denote  the true distribution of the data.

In general, model selection  can be approached from two main philosophical
perspectives: (1) testing-based model assessment and (2) parsimony-based model
assessment. The
testing-based approach involves testing, for each $k$,   $H_0: F_\tau \in \mathcal M_k$ vs $H_1:
  F_\tau  \in \mathcal M_{k+1}$.
 The Likelihood Ratio Test (LRT) is often used to build a testing-based model
 selection method. However, we do not favor a simple testing approach,
 because, in the spirit of George Box, our
 prior belief is that every restricted model is flawed, in which case every
 model would be rejected by testing at some sample size. Under this belief, a
 suitable model selection theory should  instead be based on the adequacy
 of the model approximation.

In parsimony-based model selection, we define a selection criterion, 
which usually involves a term accounting for the goodness of fit and  a
term  penalizing richer  models (models with more parameters). Then we choose as the best model
 among  a given subset of models, i.e. $\mathcal M_{best} \in  \mathbb{M}  $,  the
 model which optimizes
 our chosen criterion.  The information criterion based approaches such as AIC
 \citep*{akai:1973}, BIC \citep*{schw:1978}   and
 other model selection criteria based on approximations of Bayes factors \citep*{haug:1988} all fall under
 this large class of parsimony-based model selection tools. 
These tools are usually designed under a specific probability structure  or
are designed to achieve specific goals. For example AIC is designed to achieve
maximum predictive accuracy whereas the BIC criterion evaluates the posterior
probability of the competing models with specified priors.\citep{yang:2005} 

In practice the most parsimonious model is determined by the model associated with
minimum risk under the specific loss function. 
For example, in case of the AIC criterion, the loss function of the risk is the
Kullback-Leibler (KL) distance. We will here consider model selection from
the same basic framework but we will replace Kullback-Leibler distance with
quadratic distance \citep{Lind:2006}. We do so because the family of quadratic distances
  allows for a very straightforward
estimation of risk when the distance itself is used as the loss
function. As  evident later in this paper, the
resulting quadratic risk has many desirable properties especially for high-dimensional
data and {for} model comparisons in a non-regular parametric setup.

Before describing the quadratic risk we want to point out two underlying
themes of this article. First, as  our risk measure is in the spirit of AIC,
we will draw parallel to AIC and wherever necessary, contrast our quadratic risk
measures to AIC. Also, this research was originally motivated by the problem
of model selection for high-dimensional data 
in the normal mixture model, a non-regular parametric class. We will use this
important model to demonstrate and address pivotal
mathematical and computational problems  {in model selection}.

\subsection{ Quadratic risk based model selection}
\label{sec:quadrisk}
We will develop our
model selection criteria within the following risk framework. The loss incurred
in using a model element $M_{\theta }$ (from a model $\mathcal M $) when the
true distribution is $F_{\tau }$ will be denoted by $a(F_{\tau },M_{\theta }).$
The function $a$ should measure in some meaningful fashion the price of
using $M_{\theta }$ to approximate $F_{\tau }.$ The risk of using the model $%
\mathcal{M},$ given a particular parameter estimator $\hat{\theta},$ will be
defined to be%
$$\rho_a(F_\tau,\mathcal{M},m)=E (a(F_\tau,M_{\hat{\theta}})),$$
where the $m$ on the left side indicates that the estimate $\hat{\theta}$  on
the right side was based on a sample of size $m$,  and the expectation is
  taken over the true distribution. We alert the reader that we are retaining $m$ as a risk
argument because we will consider the estimation of risk for values of $m$
other than the actual data size $n$. 
Our goal here is to choose the model that minimizes this risk, noting that the lowest
risk model is likely to depend on $m$.  
That is, in this formal definition the problem is not to pick a true model, or
even a closest-to-true model, from $\mathbb{M}$: it is to pick a model for which
our estimated model element $M_{\hat{\theta}}$ will be the most similar to $F_\tau$, on
average, in samples of size $m$. Although models with many parameters are likely
to be closer to the true model, the estimation of the model parameters creates
extra variability that is penalized in the risk calculation.

The AIC model selection criterion  arises in this context
as follows.
If we assume both $M_\theta$ and $F_\tau$ to be continuous density
functions, with densities $m_\theta(x)$ and $f_\tau(x)$ respectively, the
Kullback Leibler(KL) loss function is
$$\mbox{KL}(F_{\tau },M_{\theta }) =\int \log\left(\frac
  {f_\tau(x)}{m_\theta(x)}\right)f_\tau(x) dx .$$
Using twice the Kullback-Leibler as a loss function, we can  split the risk into two pieces:
\begin{eqnarray}
\rho_{\mbox{\tiny KL}}(\mathcal M,m) &=& -2 E \left[\int \log(m_{\hat{\theta}}(x))
f_\tau(x) dx + \int \log(f_\tau(x)) f_\tau(x) dx \right] \nonumber\\
 &=& -2 E \left[\int \log(m_{\hat{\theta}}(x))
f_\tau(x) dx\right] + 2\int \log(f_\tau(x)) f_\tau(x) dx 
\label{eq:aic1}  
\end{eqnarray}
Here again $m$ represents the sample size used for the estimator
$\hat{\theta}$  and we will think of it as a variable that changes the magnitude of the penalty function. Of course, the convention is to replace $m$ with $n$, the actual sample size. The AIC
criterion is based on estimation of the first term on the right in \eqref{eq:aic1}. This can be done because the second term does not
depend on any of the models under consideration. As a consequence, one does
not estimate the full risk of using a model, but only its relative risk,
which in turn depends on the other models under consideration.
That is, a model that is chosen as the best among a class of models will still have an uncertain amount of total risk, and so could in fact be a poor fit.

Our quadratic risk is defined by  using  a quadratic distance  \citep{Lind:2006} as the loss function $a(.,.)$. A quadratic distance between the true distribution $F_\tau$ and any  proposed distribution $G$, where $G$ can be either discrete or continuous, has the form 
\begin{equation}
d_K(F_\tau,G)=\iint K_G(s,t)d(F_\tau-G)({s})d (F_\tau-G)(t),
\end{equation}
where $K_G$ is a kernel of appropriate dimensions (details in Section 2), that
can be chosen by the user to meet specific goals. Note that we will allow the kernel
to depend on $G$, the second argument in distance $d_K(F_\tau,G)$. 

 A familiar example of this
type is the Pearson Chi-squared distance on a partitioned sample space.
\begin{example}\label{eg:pchisq}
Let $A_{1},...,A_{C}$ be a partition of the sample space $\mathcal{S}$ into $%
C $ bins and {let $G$ be a target probability measure}. Define the kernel by 
\begin{equation}\label{eq:pcsdist}
K_G(x,y)=\sum_{i=1}^{C}\frac{\mathbb{I}[x\in A_{i}]\mathbb{I}[y\in A_{i}]}{%
G(A_{i})}, 
\end{equation}%
where $\mathbb{I}$ is the indicator function and $G(A_i)=\int_{A_{i}} dG(x)$
. The resulting quadratic
distance is the Pearson Chi-square, given by:
\begin{equation*}
\sum_{i=1}^{C}\frac{\left[ F_\tau(A_{i})-G(A_{i})\right] ^{2}}{G(A_{i})}.
\end{equation*}
\end{example}

We will be using the Pearson Chi-squared example  throughout the paper as a way to exemplify our distance calculations in the simplest possible setting. In the example we introduce in Section~\ref{selection} we will be using a Gaussian kernel.


\subsection{ Why  quadratic risk based model selection?}

Why do we propose this new measures of model selection?
Here we discuss three key issues that played a role in the
development of our methodology:

\subsubsection{ Local vs global comparison of models}
 Quadratic distance does not separate into two parts in the manner that
 Kullback Leibler does (see \eqref{eq:aic1}). This means that one must
 estimate the absolute risk for a model.  In the process, however, we can
 learn whether the models provide a low level of absolute risk, whereas for AIC we have no assurance that any model, even the best one, fits well in an absolute sense. 

The absolute quadratic risk of a parametric model can be calibrated by
comparing its risk with the risk of the empirical distribution function. The
latter estimator has no lack of fit but, due to its flexibility, but will have the highest variability. As we will see in our simulations, its risk provides a excellent benchmark for model quality. {In effect, failure to meet this standard means that either more model building is needed, or nonparametric approaches should be used. In our application section we provide
examples where the global comparison provides crucial decision making information in model selection.}





 \subsubsection{The role of model regularity conditions}
 
 The approximations used to arrive at  AIC, BIC and most other information criterion  depend
 strongly on  regularity conditions (boundary conditions, nested parameter
 structure among 
 others) some of which are violated in important model selection exercises. 
 For example, the regularity problems for nested  mixture models are well known 
(see Section \ref{selection}).  Similar problems occur in the structural equations model
 \citep{boll:2006} and the multilevel model, where the boundary between nested models is irregular. All these makes the  usual asymptotic
 expansions inappropriate. 
 
 On the other hand, our  method is based on global
 tests of goodness of fit. As a result the test statistics that are used  have asymptotic expansions
 that do not depend on regularity assumptions for nested models and so have wider validity.

 \subsubsection{ Addressing computational challenges in high dimensions}
 A primary goal of our research was to devise  model 
 selection methods that would be practical for high dimensional problems.
 However, quadratic distances, such as the one we are proposing as the loss
 function of our risk, could require 
numerical  integration of the same dimension as the data vector, which would defeat our purpose. 
For our mixture example we will show how to avoid 
this
 computational burden by using a ``rational'' kernel for the distance;
 this gives the needed integration in a closed form. A detailed
 description of construction of these kernels can be found in
 \citet{Lind:2006}.

An additional challenge in high dimensions is the construction of distances
whose operating characteristics can be tuned to the dimension of the problem
and the sample size. 
In a Chi-squared test this would be done by choosing the number and location
of the bins. 
\citet{Lind:2006} defined the degrees of freedom of a quadratic distance and discussed its use to select distances. 
The kernel we will use in the mixture example includes a
 ``tuning parameter'', that allows us to select a suitable degrees of freedom
 for the distance and also allows us to analyze the data at multiple resolutions.

\subsection{ Description of Paper}

Section~\ref{sec:distance} provides a  detailed description of our choice of
the quadratic distance in defining the quadratic risk, the  appropriateness
of such a definition in the context of high
dimensional datasets and the estimation of quadratic risk.
In Section~\ref{sec:risk-model} we build model selection tools based on our
quadratic risk. Next, in Section~\ref{selection}, we demonstrate how this model
selection tool can be applied 
to select the number of components in a multivariate normal mixture.  Application to real data and simulated datasets and comparison of
our method with  existing model selection tools are described in Section~\ref{sec:simul}.
Section~\ref{sec:discussion} contains a discussion
and the Appendix provides proofs of results stated in the article.

\section{Quadratic Risk and estimation of Quadratic Risk}\label{sec:distance}
This section provides the foundational details and the estimation techniques
for quadratic risk. 
We will start by defining the quadratic distance between two  arbitrary probability
measures.  Then,  after deriving an unbiased estimator of the quadratic
distance we will use it to build an unbiased
cross-validation based estimator for the quadratic risk
measure. Although 
analytically attractive, this cross-validation method is computationally very
expensive whenever parameter estimation is expensive, especially so for
large datasets  or when the dimensions are high. So, in the later part of this section we will derive an
AIC-like approximation to the the risk measure, which is largely based on the decomposition
of our quadratic risk and the asymptotics of quadratic distance.
The
definitions and results on quadratic distance that are discussed in this
section are described in greater detail  in \citet{Lind:2006} and \citet{Ray:2003}.


\subsection{Quadratic Distance: Definition and empirical estimate}
\subsubsection{Quadratic Distance framework}
 First we provide the details of  construction of the class of quadratic
 distance, which we use as the loss function of our risk measures.
 We define $d_K(F_\tau,G)$ to be the distance between the true density, $F_\tau$
 and any proposed distribution $G$,  provided that they are defined on the same sample space
$\mathcal{S}$.  The building block for our
distance will be $K(s,t)$, a bounded symmetric kernel function on $\mathcal{S%
}\times \mathcal{S}$, which is conditionally non-negative (CNND)\footnote{
\vskip -2ex 
\hrule
\vskip 1ex 
~~K is CNND if $\iint
K(s,t)$ $ \d \sigma (s) \d \sigma (t)\ $is nonnegative for all bounded signed
measures $\sigma ,\ $ and if
nonnegativity holds for all $\sigma $ satisfying the condition $\int  \d \sigma
(s)=0.$ }
definite. 

\begin{definition} \label{def:quaddist}
Given a CNND
$K_{G}(s,t),$ possibly depending on $G,$ the $K$-based
quadratic distance between two probability measures $F_\tau$ and $G$ is defined
as 
\begin{equation} \label{eq:quaddist}
d_{K}(F_\tau,G)=\iint K_{G}({s},t)\ \d (F_\tau-G)({s})\d (F_\tau-G)(t).
\end{equation}
\end{definition}
Note that by allowing $K$ to depend on $G$  we no longer have an inner product
space, but we retain non-negativity due to CNND.


\subsubsection{Estimation of Quadratic Distance}

We now focus on the empirical estimation of the quadratic
distance. 
 We will later use these results for the the estimation of quadratic risk.
It will be useful to express our results in terms of the $n \times n$ dimensional empirical
kernel matrix $\mathbb K$  of a data set $x_1,\ldots,x_n$, which we define to have $ij^{th}$ element
$K_{ij}=K(x_i,x_j)$.
Crucial to the estimation is the concept of the centered kernel defined as follows:
\begin{definition}
The $G$-centered kernel $K$, denoted by $K_{cen(G)}$ is defined
as 
\begin{equation}
K{_{cen(G)}(x,y)=K(x,y)\!-K(x,G)\!-K(G,y)+K(G,G)\label{eqn:ktilde}},
\end{equation}
where $ K(x,G)=\int K(x,y)dG(y)$ and $ K(G,G)=\iint K(x,y)dG(x)dG(y).$
\end{definition}
For example, the $G$-centered kernel for the Pearson Chi-squared kernel
in~\eqref{eq:pcsdist} simplifies to
\begin{eqnarray}
  \label{eq:pcsdist_cen}
 K_{cen(G)}(x,y)=\sum_{i=1}^{C}\frac{\mathbb{I}[x\in A_{i}]\mathbb{I}[y\in
     A_{i}]}{G(A_{i})} -1. 
\end{eqnarray}
Using the $G$-centered kernel we can rewrite the distance in form of a
$U$-functional  as
\begin{equation}
d_{K}(F_\tau,G)=\iint K_{cen(G)}(x,y)\d F_\tau(x)\d F_\tau(y). \label{eq-ufunctional}
\end{equation}
{But as $F_\tau$ is unknown we will use the empirical cdf $\hat{F}$ to estimate $d_{K}(F_\tau,G)$.}
 Using \eqref{eq-ufunctional}, for any fixed $G$, we  arrive at  the following two estimates for
$d_{K}(F_\tau,G)$. First, $d_{K}(\hat{F},G)=K_{cen(G)}(\hat{F},\hat{F}):=V_{n}$ is
a $V$-statistic \citep{Serf:1980} which is known to provide a biased estimate for the   $U$-functional in \eqref{eq-ufunctional}.  It can be calculated in matrix form as $\mathbf{1}^{T}%
\mathbb{K}_{cen(G)}\mathbf{1}/n^{2}.$ One can also construct an unbiased estimate  using
the corresponding $U$-statistic: 
\begin{equation}
U_{n}=\frac{1}{n(n-1)}\sum_{i}\sum_{j\not=i}K_{cen(G)}(x_{i},x_{j})=\frac{1}{n(n-1)}\left[\mathbf{1}^{T}\mathbb{K}_{cen(G)}\mathbf{1}-tr(\mathbb{K}_{cen(G)})\right],
\label{eq:unbiased}
\end{equation}%
where $tr$ denotes the usual matrix trace. Note that we will later use the notation $trace_G(K)$ to refer to a functional version of the trace operation that is defined by
 \begin{eqnarray}
trace_G(K)=\int K(x,x) \d G(x). \label{eq:theo-trace}
\end{eqnarray}
Note that for the Pearson Chi-squared distance, using the estimator $V_n$ we get
\begin{equation}
\sum_{i=1}^{C}{\left[ \hat{F}(A_{i})-G(A_{i})\right] ^{2}}\Big/{G(A_{i})},\label{eq:pcs-biased}
\end{equation}
which equals the Pearson Chi-squared statistic divided by $n$,
whereas the unbiased estimator using \eqref{eq:unbiased} gives \vskip -2ex ~
\begin{equation}
\sum_{i=1}^{C}{\left[ \hat{F}(A_{i})-G(A_{i})\right]
  ^{2}}\Big/{G(A_{i})}-\frac{(C-1)}{n}.\label{eq:pcs-unbiased}
\end{equation}
{Note that under the  usual assumptions, $n$ times
  \eqref{eq:pcs-biased} follows a $\chi^2_{C-1}$ distribution, so its expectation is  $C-1$,
  whereas the unbiased estimator in  \eqref{eq:pcs-unbiased} has expectation zero.
}

\subsection{ Quadratic Risk definition and unbiased estimation} \label{sec:exact-risk}
Based on the quadratic distance we have already defined the quadratic risk of a parametric model $\mathcal M$  as 
\begin{equation}
  \label{eq:quadrisk}
 \rho_d(F_\tau,\mathcal{M},m)=E (d_K(F_\tau,M_{\hat{\theta}})).
\end{equation}
Here $\hat{\theta}$ is an estimator of ${\theta}$, based on a sample of size
$m$. In our subsequent analyses
we will assume that  $\hat{\theta}$ is the maximum likelihood estimator of ${\theta}$.
Note also that  $M_{\hat{\theta}}$ is the second argument in $d_K(.,.)$, so when we
use a $G$-dependent kernel (see \eqref{eq:quaddist}) it will be playing the
role of $G$. 


Next we construct an unbiased estimator of the risk  given in
\eqref{eq:quadrisk}. 
We start by showing how
to unbiasedly estimate the quadratic risk of the empirical distribution
function, a nonparametric estimator of $F_\tau$. 
As mentioned earlier, comparison of nonparametric risk and parametric risk gives a global assessment of model quality.

\subsubsection{Unbiased estimation of empirical risk}\label{sec:empirical}
For simplicity of notation, we define the risk of the nonparametric fit,  based on $m$ observations, by
$R_m$. A straightforward calculation shows that
\begin{eqnarray}
R_m = E [ d_K(F_\tau,\hat{F})]=\frac 1 {m} trace_{F_\tau}(K_{cen(F_\tau)}).
\label{eq:Rn}
\end{eqnarray}
For example, for the Pearson Chi-squared example the above risk can be calculated as:
$$ \frac 1 {m}      \sum_{i=1}^C \left[ \frac 1 {F_\tau(A_i)} - 1
  \right] F_\tau(A_i) = \frac 1  m \sum_{i=1}^C [ 1- F_\tau(A_i)] =  \frac{C-1}{m} .
$$
For most other kernels we will need to estimate $R_m$ because it will depend on $F_\tau$.
Provided that $K$ does not depend on $G$, we can compute an unbiased estimator of
$trace(K_{cen(F_\tau)})$, based on a sample of size $n$  using
\begin{eqnarray}
\label{eq:rn-1}
 \frac 1 n \sum_{i=1}^n {K}(x_i,x_i)- \frac 1 {n(n-1)}
\sum_i\sum_{j\not=i} {K}(x_i,x_j).
\end{eqnarray}
This immediately gives us the unbiased estimator of $R_m$ as 
\begin{eqnarray}
\label{eq:rn-2}
\hat{R}_{m}&=& \frac 1 {m} \Big[ \frac 1 n \sum_{i=1}^n {K}(x_i,x_i)- \frac 1 {n(n-1)}
\sum_i\sum_{j\not=i} {K}(x_i,x_j)\Big].
\end{eqnarray}

Note that $\hat{R}_{m}  \xrightarrow{prob} 0$ as $n \rightarrow
  \infty$
and if we set $m=n$, then $\widehat{R_n}=O_p(1/n)$. Hence any method of systematically selecting models with smaller estimated risk than $\widehat{R_n}$ will give risk estimates converging to zero.

\subsubsection{Unbiased estimation of risk of parametric models}
We next build an unbiased estimator of the risk of the parametric model
$\mathcal M= \{M_{\theta}\}$, which  is given in  \eqref{eq:quadrisk}
This risk can be estimated unbiasedly at any sample size $m\le n-2$ as
follows:
 Let $A_{m}$ be a randomly selected subset of size $m$ from $\{1,2,...,n\}$, and let the point
 estimator $\hat{\theta}(A_{m})$ be the
value of $\hat{\theta}$ based on $\{x_{i}:i \in A_{m}\}$. Also, we define 
\begin{equation}
U(A_{m})=\frac{1}{(n-m)(n-m-1)}\sum_{{ i,j\in A_{m}^\complement, i\neq j}}
K_{cen\left(M_{\hat{\theta}(A_{m})}\right)}(x_{i},x_{j}),
\label{eq:uam}
\end{equation}
$A^\complement$ being the complement of set $A$.
As $\hat{\theta}$ is constructed from $m$ independent observations, $U(A_{m})$ is an unbiased estimator of $\rho (\mathcal M, m),$ for $m\leq n-2.$
Using \eqref{eq:uam} we can define an unbiased  estimator   by
averaging over all possible 
subsets of size $m$,
$$ \hat{\rho}(F_\tau,\mathcal M,m)= \frac{1}{\left(\substack{ n \\ m} \right)}
\sum_{l=1}^{\left(\substack{ n \\ m} \right)} U(A_{{m}_l}). $$

If one wishes  to estimate $\rho(F,M,n)$, then  $m=n-2$ might be expected to generate the least bias.
 However, as we argue at the end of this section, other 
choices of $m$ might be made for the purposes of either consistency or of parsimony. 
Also, for computational reasons one might wish to use only a few selected subsets of size $m$, as 
in V-fold cross-validation  \citep{vander:2004}, or randomly selected subsets of size $m$ \citep{blom:1976}.


Note that one can construct unbiased estimators of the AIC relative risk
 in a similar fashion. That is, one can estimate   the first term on the
 right hand side of \eqref{eq:aic1} as $ \sum \log \left[m(\hat{\theta}_{(-i)})(x_i)\right]$. This estimator is
 unbiased for risk at sample size $m=n-1$.

\subsection{ Decomposition  of quadratic risk and approximate quadratic
  risk estimators}
\label{sec:decomp}

{In many model selection problem, we may choose to} simplify the unbiased estimates of
risk, which can be both difficult and expensive to compute. The main difficulty
is that the point estimators are traditionally obtained using the
notoriously slow EM algorithm 
and therefore it is difficult to obtain them to sufficiently high precision
in each of the many ``leave-out'' recalculations.  
For this reason, we derive an
alternative estimate of the quadratic risk, based on asymptotic expansions,
similar to those of the AIC derivation.

We start with the following decomposition of  the
distance:
\begin{eqnarray}
    d_K(F_\tau,M_{\hat{\theta}}) & = 
& d_K(F_\tau,M_{{\theta}_\tau})+\left[d_K(F_\tau,M_{\hat{\theta}})-
  d_K(F_\tau,M_{{\theta}_\tau})\right]\label{eq:distdecomp}.
 \end{eqnarray}
Here ${{\theta}_\tau}$ denotes the parameters of the distribution (in the class
of distributions denoted by $M_\theta$) closest in quadratic
distance to  the true
model $F_\tau$. That is, $\theta_\tau\!\!=\!\! \argmin_\theta   d_K(F_\tau,M_{{\theta}})$.
We will call the first term on the right side of \eqref{eq:distdecomp},
the model lack of fit (MLF)
because it assesses the distance between the best model element and the
truth. In particular, it is  zero  if $F_\tau \in M_\theta$. 
Notice that MLF does not depend on parameter estimation.

We will call the  second term  on the right side of \eqref{eq:distdecomp}, 
$\left[d_K(F_\tau,M_{\hat{\theta}})- d_K(F_\tau,M_{{\theta}_\tau})\right]$,
the parameter estimation error because it measures the deviation of $\hat{\theta}$ from the best parameter $\theta_\tau$.  
We observe that it is always non-negative, by definition of $\theta_\tau$,
 and its magnitude increases as
$\hat{\theta}$ deviates from $\theta_\tau$. 
{In fact, it is  shown in the appendix
(equation \eqref{eq:quadloss}) that this term is approximately a quadratic form in
$(\hat{\theta}-\theta_\tau)$, which implies that the parameter estimation cost is a monotonic function of $(\hat{\theta}-\theta_\tau)$}. Based on the interpretation of the decomposition of the distance
we can rewrite \eqref{eq:distdecomp} as
\begin{eqnarray}
    d_K(F_\tau,M_{\hat{\theta}}) & = 
& \mbox{MLF}+\mbox{ Parameter Estimation Error}.\label{eq:distdecomp1}
 \end{eqnarray}
For the risk calculation we take expectations in \eqref{eq:distdecomp1} to get
\begin{equation}
 E(d_K(F_\tau,M_{\hat{\theta}})) = \mbox{MLF}+\mbox{PEC}_{(m)},\label{eq:riskdecomp}
\end{equation}
where $PEC_{(m)}=E\left[d_K(F_\tau,M_{\hat{\theta}})-d_K(F_\tau,M_{{\theta}_\tau})\right]$  represents the parameter estimation
cost when ${\hat{\theta}}$ comes from a sample of size $m$.

This decomposition is essential both for finding an approximation to the
overall risk and for providing a useful interpretation of the different errors
driving the risk. 
The approximation and the estimates for the two terms, MLF and
PEC, will be directly based on the
asymptotic theory of quadratic distance \citep{Lind:2006} and the unbiased
estimators of quadratic distance.
 
Now we will  state a few results from \citet{Lind:2006} and  use these results to justify our approximation and finally provide
an appropriate estimator for the quadratic risk.

\subsubsection{ Kernel projection operator and score centered kernel}
Crucial to these approximations of quadratic risk are the concept of 
scored centered kernel and the score based projection operators which we now define.

For a parametric model given by $G_\theta$, with density given by
$g(x;\theta)$, let us denote the set of score
functions by  $\nabla \log g(x_i;\theta) = \mathbf{s}(x_{i};\theta ),$
where $\nabla$ is the vector differential with respect to every element of  $\theta$.
Note the MLE $\hat{\theta}$ is obtained as a solution to $\sum
\mathbf{s}(x_{i};\theta )=0$. Further, we define the extended score
vector {$\mathbf{u}^{\ast }_\theta(x)=(1,\mathbf{s}(x,\theta)^T)^T$} and the extended information matrix for a single observation to be: 
\begin{equation*}
\mathbb{J}_{\theta }^{\ast }=E_{\theta }[\mathbf{u}_{\theta }^{\ast }\mathbf{%
u}_{\theta }^{\ast T}].
\end{equation*}
We will then let $P^{\ast }$ be the kernel operator defined by%
\begin{equation}
P_{\theta }^{\ast }(x,y)=\mathbf{u}_{\theta }^{\ast }(x)^{T}\cdot \mathbb{J}%
_{\theta }^{\ast -1}\cdot \mathbf{u}_{\theta }^{\ast }(y). \label{eq:P=uju}
\end{equation}%
%
That is, $ P_{\theta }^{\ast }$ acts as a projection operator on the
likelihood scores $\mathbf{u}_{\theta }^{\ast}$.
The score centered kernel $K_{scen(G_\theta) }(x,y),$ as centered under $%
G_{\theta },$ is defined to be 
\begin{eqnarray}
K_{scen(G_{\theta })} &= &(I-P_{\theta }^{\ast })K(I-P_{\theta }^{\ast })
\nonumber \\
&=& 
K(x,y)-\!\!\int\! \!P_{\theta }^{\ast }\mathbb{(}x,z)K\mathbb{(}z,y)\ dG_{\theta
}(z)
-\!\!\int \!\! K(x,z)P_{\theta }^{\ast }(z,y)dG_{\theta }(z)
\nonumber \\
&&\qquad+\!\iint \!\!P_{\theta
}^{\ast }(x,z)K(z,w)P_{\theta }^{\ast }(w,y)dG_{\theta }(z)dG_{\theta }(w).
\label{score centering}
\end{eqnarray}

Similar to the matrix form of $\mathbb K$, we define the 
$n\times n$ matrix $\mathbb{%
P}_{\theta }=$ $\mathbf{u}_{\theta }\mathbb{J}_{\theta }^{-1}\mathbf{u}%
_{\theta }^{T},$ where \textbf{u}$_{\theta }$ is the $n\times p$ matrix with
entries $\partial _{\theta _{j}}[\log g_{\theta }(x_i)], x_i$ being the
$i^{th}$ data point and $\theta _{j}$ being the $j^{th}$ component of $\theta$. Similarly we
represent the matrix version of the scored centered kernel as
 $\mathbb{K}_{scen({\theta })}$ which can be calculated as
$$   \mathbb K_{scen({G_\theta })} = (\mathbb{I}-\mathbb{P}^{\ast}_{\theta
})\cdot \mathbb{K}\cdot (\mathbb{I}-\mathbb{P}^{\ast}_{\theta }) 
= (\mathbb{I}-\mathbb{P}_{\theta
})\cdot \mathbb{K}_{cen({G_\theta })}\cdot (\mathbb{I}-\mathbb{P}_{\theta }).
$$

The first important property of score centering is the following alternative
representation for the empirical distance between the data and the estimated
model:%
\begin{equation*}
d_{K}(\hat{F},G_{\hat{\theta}})=\iint K_{scen{(G_{\hat{\theta}})}}(x,y)\  \d \hat{F
}(x) \d \hat{F}(y),
\end{equation*}
which has the following asymptotic property:
\begin{theorem}
\label{thm:convergence} \textbf{\citet{Lind:2006}}
Given the regularity assumptions stated in \citet{Lind:2006}, under $G_{\theta }$
we have 
\begin{equation}
n\left[ d_{K}(\hat{F},G_{\hat{\theta}})-\iint K_{scen(G_\theta) }(x,y) \d\hat{F}(x) \d \hat{F}(y)\right] \xrightarrow{prob} 0.
\end{equation}%
\end{theorem}
This result can then be used to show that $nd_K(\hat{F},G_{\hat{\theta}})$
has an asymptotic distribution that can be represented as $\Sigma_{i=1}^\infty \lambda_i
\chi_{1}^2$, where the $\lambda_i$'s are the eigenvalues of the spectral
decomposition of kernel $K_{scen(G_\theta)}$. The asymptotic mean is therefore
$\Sigma_{i=1}^\infty \lambda_i=trace \left(K_{scen(G_\theta)}\right)$ \citep[see][for details]{Lind:2006}.

\subsubsection{Approximate quadratic risk estimators}
We now return  to our original purpose of approximating the different
terms of our quadratic risk. These approximations will be based on asymptotic
calculations for each model,
 assuming that the model is correct. This  same technique is used in the AIC calculation, where it clearly produces great simplification in risk estimation. We
start with the  asymptotic approximation to PEC: 
\begin{PRO}
 \label{pro:pec}
\label{pro:pec} 
Suppose  that $M_{\theta }$ has the density function $m_{\theta }$ and that
the true model $F_{\tau }$ is in the model space $M_{\theta },$ having
parameter $\theta _{\tau }.$ In this case $d_{K}(F_{\tau },M_{{\theta }%
_{\tau }})=0.$ Define 
\begin{equation*}
V(x,y)=\iint P(x,z)K(z,w)P(w,y)\ M_{\theta _{\tau }}(z)M_{\theta _{\tau }}(w)
\end{equation*}%
where $P(x,y)=\mathbf{u}(x)^{T}\mathbb{J}^{-1}\mathbf{u}(y).$ Then as $m,$
the sample size for $\hat{\theta},$ goes to infinity,   
\begin{equation*}
md_{K}(F_{\tau },M_{\hat{\theta}})-m\iint V(x,y)\ \hat{F}(x)\hat{F}%
(y)\xrightarrow{prob}0.
\end{equation*}%
It follows that the asymptotic mean of $md_{K}(F_{\tau },M_{\hat{\theta}})$
is 
\begin{equation*}
tr_{\theta _{\tau }}(V)=\iint V(x,y)\ M_{\theta _{\tau }}(x)M_{\theta _{\tau
}}(y).
\end{equation*}
\end{PRO}
\begin{proof}
See Appendix.
\end{proof}

\noindent
As a simple example, for the Pearson Chi-squared kernel the PEC
is simply
$dim(\hat \theta)/m$, where $dim(\hat \theta)$ denotes the number of parameters being
estimated.
In general we have to estimate PEC and for an estimator  of  $trace(P^\ast_{\theta } K P^\ast_{\theta})$, we use  its empirical 
 version  {$\frac 1 n tr (\mathbb{P^\ast_{\hat{\theta}} KP^\ast_{\hat{\theta}}})
$ based on the full sample. This yields the following estimator: 
\begin{equation}
\widehat{\mbox{PEC}}_{(m)}= \frac 1 m \left( \frac 1 n tr (\mathbb{P^\ast_{\hat{\theta}} KP^\ast_{\hat{\theta}}}) \right)\label{eq:pecest}.
\end{equation}

Now we turn  to the estimation of MLF=$d_K(F_\tau,M_{{\theta}_\tau})$,  a term not depending on hypothetical sample size $m$. We start with its sample equivalent given by
\begin{equation}
  \label{eq:mlf1}
 d_K(\hat{F},M_{\hat{\theta}})=\iint K_{scen(M_{\hat{\theta}})}(x,y)  \d \hat{F%
}(x) \d  \hat{F}(y)  = \frac 1 {n^2} { \mathbf{1}^T\mathbb{K}_{scen(M_{\hat{\theta}})}  \mathbf{1}}.
\end{equation}
This estimator has bias  $E(d_{K}(\hat{F},M_{\hat{\theta}}))-MLF$. To do bias
correction we again use an
asymptotic approximation assuming that $F_\tau=M_{{\theta}_\tau}$, so that
$MLF=0$. Applying the remarks following Theorem~\ref{thm:convergence} we obtain the following
approximation to the bias term:
\begin{equation}
  \label{eq:mlf2}
 E(d_K(\hat{F},M_{\hat{\theta}}))-MLF\approx
\frac 1 n trace_{\theta_\tau}(K_{scen({M_{\theta_\tau}})}).
\end{equation}
The right side of   \eqref{eq:mlf2} can be estimated by its sample equivalent
\begin{equation}
\frac 1 {n^2}
tr(\mathbb{I}-\mathbb{P^\ast_{\hat\theta}})\mathbb{K}(\mathbb{I}-\mathbb{P^\ast_{\hat\theta}}).
\label{eq:biasadjust}
\end{equation}
Thus,  we have the bias corrected estimator of MLF,
\begin{equation}
 \widehat{\mbox{MLF}}=\left[d_k(\hat{F},M_{\hat{\theta}})-\frac 1 {n^2}
  tr(\mathbb{I}-\mathbb{P^\ast_{\hat\theta}})\mathbb{K}(\mathbb{I}-\mathbb{P^\ast_{\hat\theta}})\right]. \label{eq:mlf}
\end{equation}
For example, in the Pearson Chi-squared case \eqref{eq:biasadjust} reduces to
$(C-1 - dim(\hat \theta))/n$ and thus the MLF can be estimated as:
$$\sum_{i=1}^{C}\frac{\left[ \hat{F}(A_{i})-G_{\hat{\theta}}(A_{i})\right] ^{2}}{G_{\hat{\theta}}(A_{i})} -\frac{C-1 - dim(\hat \theta)}{n}. 
$$
Since 
$(C-1-dim(\hat \theta))$ gives the residual degrees of freedom for the $C$-cell
Pearson Chi-squared test it is clear how the correction removes the bias inherent in the Chi-square limiting distribution.  

Finally combining the estimators of $\mbox{PEC}_{(m)}$ and MLF given in  \eqref{eq:pecest}
and \eqref{eq:mlf} we have the following estimator for the quadratic risk:
\begin{eqnarray}
   \widehat{\rho_{d_K(\mathcal M,m)}}&=&\left[d_k(\hat{F},M_{\hat{\theta}})-\frac 1 {n^2}
     tr(\mathbb{I}-\mathbb{P^\ast_{\hat\theta}})\mathbb{K}(\mathbb{I}-\mathbb{P^\ast_{\hat\theta}})\right]+\left[\frac 1 {nm} tr (\mathbb{P^\ast_{\hat\theta}}\mathbb{K}\mathbb{P^\ast_{\hat\theta}})\right]\nonumber\\
&=&\left[
\frac 1 {n^2}
\mathbf{1}^T\mathbb{K}_{scen(M_{\hat{\theta}})}  \mathbf{1}
-\frac 1 {n^2}
     tr(\mathbb{I}-\mathbb{P^\ast_{\hat\theta}})\mathbb{K}(\mathbb{I}-\mathbb{P^\ast_{\hat\theta}})\right]+\left[\frac 1 {nm} tr (\mathbb{P^\ast_{\hat\theta}}\mathbb{K}\mathbb{P^\ast_{\hat\theta}})\right]
\label{eq:riskfinal}
\end{eqnarray}

Note that  for the Pearson Chi-squared example the risk estimator, when
calculated at $m=n$  simply becomes:
\begin{equation}
\sum_{i=1}^{C}{\left[ \hat{F}(A_{i})-G_{\hat{\theta}}(A_{i})\right]
  ^{2}}\Big/{G_{\hat{\theta}}(A_{i})}-\frac{(C-1)+2 dim(\hat \theta)}{n}. \label{eq:PCSrisk}
\end{equation}

Observe that   the risk estimator given by \eqref{eq:riskfinal} is
numerically less expensive than the cross validation based exact estimator given in
\eqref{eq:uam}.
 Unlike the unbiased risk estimator, it does not require repeated parameter estimation.
Moreover, the matrix operations involved in the calculation  are
also inexpensive, making  \eqref{eq:riskfinal} a promising risk estimator.

For model comparison purposes we will also want to estimate the risk of the empirical
distribution function $\hat{F}$, in which case  $\mathbb{P^\ast_{\hat\theta}}=I$. The estimated  distance $d_K(\hat{F},\hat{F})$ is zero. 
So a biased 
estimator of $R_m$, the risk of the  $\hat{F}$, at  sample size $m$, may
be calculated using
\begin{equation}
  \label{eq:Rn-biased}
\widehat{R^b_m}= \frac 1 {nm} tr \left( \mathbb K_{cen(\hat{F})}\right) .  
\end{equation}

{We will later see that if we use an appropriately scaled kernel $K^s$,
  \footnote{ 
\vskip -2ex 
\hrule
\vskip 1ex 
~~$K^s=K/c$, where $c=tr(K)/tr(K^2)$   } $\widehat{R^b_m}$ is exactly equal
to the degrees of freedom corresponding to the scaled kernel $K^s$. In fact for
the Pearson
Chi-squared kernel $\widehat{R^b_m}=C-1$, which is the degrees of freedom
for a $C$-cell Chi-squared goodness of fit test.
}

In the special case that the estimated risk is evaluated at $m=n$ we will
call it the Quadratic AIC risk or in short QAIC, where AIC reflects commonality with the AIC
derivation. We will now show how to  mimic the BIC criterion by using a different value of $ m$. 
 Recall that for the usual AIC
with KL loss function,  the  relative MLF, after subtracting $\int \log(f_\tau(x)) f_\tau(x) dx$ common to all models  corresponds to $-\!\! \int \log(m_{\theta_\tau}) f_\tau(x)  $ 
(see \eqref{eq:aic1}), which has the approximately unbiased estimator (under the model):
\begin{equation*}
\widehat{\mbox{MLF}}= -2\frac{\hat{l}} n + \frac{ dim(\hat \theta)} {n} .
\end{equation*}
The parameter estimation cost, $\widehat{\mbox{PEC}_{(m)}}= { dim(\hat \theta)}/{m}  \label{eq:pecaic},$ giving us the estimated risk,  
\begin{equation}
-2\frac{\hat {l}} n
  + \frac{ dim(\hat \theta)} {n} + \frac{ dim(\hat \theta)} {m} .\label{eq:aicdeco}
\end{equation}
{For AIC risk we use $\mbox{PEC}_{(n)}$, which gives us $AIC/n=  -2 {\hat {l}}/ n
  + 2{ dim(\hat \theta)} /{n}$. }
Note the similarity of this risk to
the Pearson risk estimator in \eqref{eq:PCSrisk}. The essential difference is that the latter
estimates absolute risk, not relative risk. 
If one uses $m=n/(\log n-1)$ in \eqref{eq:aicdeco}, one arrives at the standard BIC formula:
$$\frac {BIC} n = -2\frac{\hat {l}} n
  + \frac{ \log(n) dim(\hat \theta)} {n} .$$ 
For this reason the estimated risk \eqref{eq:riskfinal} evaluated at
$m=n/(\log n-1)$ will be called QBIC, and can be written as
\begin{eqnarray}
  \mbox{QBIC} & = &\widehat{\mbox{MLF}}+ (\log(n)-1)
  \times\widehat{\mbox{PEC}}_{(n)}.\nonumber \\
  &=&\left[d_k(\hat{F},M_{\hat{\theta}})
-\frac 1 {n^2}
     tr(\mathbb{I}-\mathbb{P^\ast_{\hat\theta}})\mathbb{K}(\mathbb{I}-\mathbb{P^\ast_{\hat\theta}})\right]+\left[\frac  {(\log(n)-1)} {n^2} tr (\mathbb{P^\ast_{\hat\theta}}\mathbb{K}\mathbb{P^\ast_{\hat\theta}})\right]. \label{eq:qbicrisk}
\end{eqnarray}

{Along with a recipe for estimation of the risk, the decomposition 
in~\eqref{eq:riskfinal} provides an
excellent tool for deeper understanding and analysis  of the interplay of
model misspecification and  the effective parameter cost of using the
model. This provides us with the ingredients for  constructing alternative
model selection tools and tools for accessing global fits which we discuss in
the following section}

\subsection{Consistency and choice of $m$}\label{sec:consistency}

It seems very natural to assess the risk of a model at the sample size $m=n.$
This risk indicates how well the chosen model, along with the current
parameter estimates, might be expected to perform when used to approximate
the true distribution. Unfortunately, estimating the risk at $m=n$ is a hard
problem: one cannot estimate this risk well enough to discriminate
consistently between  two true models with different numbers of parameters. 

In this section we will consider the meaning of consistency under several
asymptotic scenarios. We will use the AIC method to illustrate our points
because of its simple form. In all of them, we will suppose that the true
distribution is in one of the models, i.e $F_\tau\in \mathcal{M}$. 

First, if one leaves $m$ fixed as $n\rightarrow \infty ,$ then the best
model, in terms of risk, need not contain the true distribution, and
consistency would mean selection of this (possibly false) best model. This
might be desirable for reasons of parsimony, as the true models may be too
complex to be useful.

In a formal sense, methods like AIC cannot be consistent in this {asymptotic scenario}
because they use asymptotic approximations for the cost of estimating
parameters, and those approximations assume the model is true. If one
replaces the exact risk with the asymptotic approximation, however,
consistency is achieved. See,  for example, the AIC in \eqref{eq:aicdeco}, where the
$\widehat{\mbox{MLF}}$ term is consistent for MLF, and the approximate parameter estimation cost is just
a function of parameter dimension when $m$ is fixed.

A second way to view consistency is to let $m=m_{n}$ increase with $n.$ In
this case, the parameter estimation costs shrink to zero, and so consistency
would involve selecting a true model. A  stringent form of consistency,
but natural, would require the method to select the smallest true model with
probability approaching one. We will call this smallest-true-model (STM)
consistency.

It is well known that AIC, with $m_{n}=n,$  lacks this form of consistency.
When two true models $\mathcal{M}_{1}\subset \mathcal{M}_{2}$ are nested (and
regularity conditions hold), the probability that AIC picks the larger of
the two models corresponds to using a likelihood ratio test with the
constant critical value $2(\dim (\mathcal{M}_{2})-\dim (\mathcal{M}_{1})).$
For example if the difference of dimension is 1, with probability
$0.16(=Pr(\chi^2_1>2))$ we choose the larger model.

It turns out that STM consistency holds when $m/n$ goes to zero. The
proportional case, when $m_{n}\approx \gamma n,$ for constant $\gamma >0,$
is at the cutting edge of consistency. It is easily checked that for the AIC
criterion, it is equivalent to using the critical value $(1+\gamma ^{-1})$ $%
(\dim (\mathcal{M}_{2})-\dim (\mathcal{M}_{1})).$ For example if the
difference of dimension is 1, the probability of overshooting is around
$.045,.025, .01,.001$ for the proportion $\gamma= \frac 1 3 ,\frac 1 4, \frac 1 5, \frac 1 {10} $ respectively.
Thus as $\gamma $ approaches zero, the probability of overshooting goes to zero.

Although one might choose to use a BIC type criterion for its familiarity,
we think it would be very useful to also describe the null probability of
rejection for the corresponding likelihood ratio critical value.

\section{Model Selection using Quadratic Risk}  \label{sec:risk-model}

\begin{figure}[tbp]
  \centering 
\includegraphics*[width=.8\textwidth]{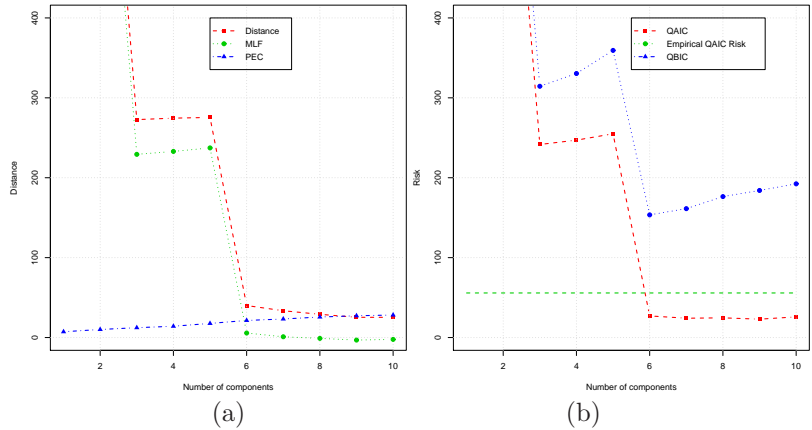}
\vskip -.2in
~\hfill\hfill(a)\hfill \hfill(b)\hfill\hfill~
\caption{Illustration  of (a) components of the risk measure and (b) and the  risk measures
    corresponding to  a single run of the
  simulation of Model(3) described in Section 5.1.1. \label{fig:3_6}} 
\end{figure}

Now, we provide a short description of how one might use quadratic risk as a model
selection tool. First, one needs to understand the interplay of the terms in
the risk decomposition.  Please see Figure~\ref{fig:3_6} for an
illustration. It shows the relevant quantities for an example in Model(3) with $n=300$, described later in Section 5. Suppose we have a sequence of models $M_k$ that are
nested with increasing $k$.
As model complexity $k$ increases the empirical distance $d_K(\hat{F},\hat{M_k})$   decreases but must stay non-negative.
We also know that $\widehat{\mbox{MLF}}$, for any model that contains $F_\tau$, will tend to zero in probability.
 The estimated parameter estimation cost $\widehat{\mbox{PEC}}$
starts near zero and increases
with model complexity. However it never will be larger than  the risk of the empirical
estimate, which is the parameter estimation cost of a non-parametric fit. By construction, the
total quadratic risk  is the simple addition of the MLF and the
 PEC. Thus models with low estimated risk arise as a compromise between the
decreasing MLF and the increasing PEC. 
A simple model selection rule is to pick the model having the minimum
value of QAIC.

The QBIC form of the risk, given in \eqref{eq:qbicrisk} introduces a larger penalty for parameter estimation
than QAIC. In Figure~\ref{fig:3_6}(b) one can see how this alteration has created a much sharper minimum in the risk function. 
 Both for QAIC and QBIC,
we will say that the model with the minimum estimated risk is the best model.

Just like AIC and BIC, quadratic risks need not attain an absolute minimum in
the range of
models considered. 
  However, there is an important benchmark, given by the empirical risk, for evaluating the performance of a model when
using  quadratic risk measures. Note  that for the empirical distribution there is no model lack of fit, because the model is nonparametric, but there is the maximal possible parameter estimation cost. 
 Models
that don't meet this benchmark clearly suffer from substantial model lack
of fit and so we will label them as ``inadequate models,'' while the rest
will be called ``adequate''.
Even in the presence of a minimum risk model, if all our proposed models are
``inadequate'' it indicates that more model building in the same class, or
exploration in a larger class of models might be necessary. For example,
 based on the QAIC risk for proposed models and empirical model in
 Figure~\ref{fig:3_6}(b), if we had explored only the 5 smallest
 models, then QAIC and QBIC would have selected  model 3. However based on
 the MRA criterion no model would have been found adequate, thus providing
 the crucial global framework for model selection.

On the other hand if we have more than one ``adequate models''
 we  will call the smallest
model in this set   the minimal-risk-adequate model
(henceforth denoted as MRA).  Note that the consistency
 of the empirical risk (shown in Section \ref{sec:empirical})ensures that this model selection method would have
 estimated risk going to 0 as $n \rightarrow \infty$
Now inspecting
 Figure~\ref{fig:3_6}(b) for the whole range of 10 components model, we see
 that all models from 6 onwards are adequate. Thus model 6 is the MRA model, even though the QAIC  selects model 9 as the minimum risk model.

  We will provide illustrations of these criteria in the application section.

\section{ Application: Selecting the number of components in
  mixture models}
\label{selection}

Now we apply our quadratic risk model selection methodology to the problem of
selecting  the number of components in a
multivariate normal mixture model.
We start by introducing the following notation and definitions for mixture models.

A random variable $ X\in\mathbb{R}^D$  is said to follow a  $k$-component
normal mixture model if its density
$f^{(k)}$ can be written as
\bea
 f^{(k)}(x;\boldsymbol\mu,\boldsymbol\Sigma,\boldsymbol\pi) = \sum_{j=1}^k \pi_j~~ \phi(x;\mu_j,\Sigma_j) \mbox { for } x\in
  \mathbb{R}^D, \label{eq:mix}
\eea
where  $\phi(x;\mu_j,\Sigma_j)$ denotes a normal density with mean $\mu_j$ and
variance $\Sigma_j$, and 
$ \pi_j < 1 \forall j$ and $\sum_{j=1}^k \pi_j = 1$. For compactness we
denote $\theta=\{\boldsymbol\mu,\boldsymbol\Sigma,\boldsymbol\pi\}$, with
the above restrictions and denote  the density in \eqref{eq:mix} as  $f_\theta^{(k)}$.

According to the general framework, now $\mathbb{M}$ denotes the class of all $D$-dimensional finite
mixtures with normal mixing components; i.e. for a fixed $D$, our
$\mathbb{M}=\{\mathcal M_{1},\mathcal M_2,\ldots\}$,
where $\mathcal M_k$  denotes the $k$-component normal mixture models,
and for a particular model element  $M_\theta \in \mathcal M_k$, more precisely denoted as
$M_\theta^{(k)}$ we have
$d M_\theta^{(k)}= f^{(k)}_\theta$.

The complexity of the null distribution makes a testing theory for the number
of components rather more difficult. It also means that the standard
regularity assumptions behind the derivations of the AIC and BIC are
false. For a description of the previous approaches to selecting the number of components of mixture models see \citet{Mcla:2000}.

\subsection{Risk-based analysis of mixture complexity}
We now develop a strategy for selecting the number of components using  our quadratic risk assessment approach. 
  First, we will
specify the kernel in order to specify the loss
function. We will then outline the steps for
estimating the risk function and finally illustrate the use
of the risk measure in choosing a mixture model.

\subsubsection{Specifying the kernel}\label{sec-kernel-specification}
In principle, any nonnegative definite kernel can be used
to build the distance measure. However, from the results in
Section~\ref{sec:distance} it is clear that the key to
the calculation of the distance, and hence the risk,
are the integrals required in the definition of the quadratic distance.
 These will involve high dimensional numerical
  integration for arbitrary choice of the kernel $K$. 
For this reason it is desirable to use kernels for which the integrations
$\int K(x,y) dM(y)$ and $\iint K(x,y) dM(x) dM(y)$
can be carried out explicitly. 

When the model is a finite mixture of normals,
the  multivariate normal kernel meets this goal. 
 The $D$-dimensional normal kernel, in its most general form, is defined as
\begin{equation*}
K_{\Sigma }(x,y)=\frac{1}{(2\pi )^{\frac{p}{2}}|\Sigma |^{\frac{1}{2}}}%
\exp\left( -\frac{1}{2}(x-y)^{\prime }\Sigma ^{-1}(x-y)\right) .\label{eq:ksigma}
\end{equation*}
We will  take $\Sigma =h^{2}I,~~h$ being a ``smoothing
parameter'' and henceforth denote the corresponding kernel as $K^{h}$. This constant $h$ can be thought of as a ``tuning'' parameter,
something analogous to the bin-width in the construction of a histogram.
We will discuss its role further in the next subsection. 

\subsubsection{ Role of tuning parameter ``$\bm{h}$'' and its empirical estimation}
\label{sec:sdof}
\citet{Lind:2006}  developed a tool for the understanding of the operating
characteristics of a quadratic distance they called the spectral degrees of
freedom of the kernel (sDOF). It is a generalization of the degrees of freedom of
the chi-squared distance (itself a quadratic distance) to other kernels by
examining their limiting distributions. 
With this tool one can roughly equate the  degree of smoothing that comes
from a choice of $h$ to that of the chi-squared test having the same degrees of
freedom, which in turn corresponds to the number of cells used in its
construction. 
The tuning parameter $h$ is
analogous to choosing the bin-width of each cell (or, equivalently  choosing the
number of cells) in a Chi-squared goodness of fit test.

The spectral degrees of freedom of the kernel depend on the true sampling distribution and the kernel $K$ (whether centered or score centered) through the formula
   \bea
{sDOF}= {\left(\int K^h(x,x) dF_\tau(x) \right)^2}\Big/ { \int \left({K^h(x,x)}\right)^2 dF_\tau(x)}.
\eea
In this paper,
we will base our selection of $h$ on the empirical estimate of $sDOF$ given by
\bea
\widehat{sDOF}= {\left(tr(\mathbb K^h_{cen(\hat F)})\right)^2} \Big/ 
{ tr \left({\mathbb K^h_{cen(\hat F)}}^2\right)} ,
\eea
where the centered kernel matrix based on the empirical distribution has the
$ij^{th}$ element
$$ \textstyle K^h_{cen(\hat F)}(x_i,x_j)= {K}_h(x_i,x_j) -\frac 1 n \sum_i
 {K}_h(x_i,x_j)-\frac 1 n \sum_i {K}_h(x_i,x_j)+\frac 1 {n^2} \sum_i \sum_j
 {K}_h(x_i,x_j). $$

 Just as in a Chi-squared goodness of fit test, there are choices for the
spectral degrees of freedom that are clearly too small and others that are clearly too
large.  As a rough  rule of thumb for selecting the
number of cells for $D=1$ we suggest  that  the degrees
of freedom should be  more than $5$ and less than $ n/5$,  $n$ being the total
number of observations. For $D>1$ dimensions one needs to increase the lower bound of 5 because the goodness of fit test now must assess fit in several directions. We find that if $sDOF$ is smaller than $\binom{D+1}{2}$, it is likely
to be  too
much smoothing, whereas if $sDOF > n/5$ it would mean too little smoothing. 
See our simulation section for more on the role of the
smoothing parameter on model selection.



\subsection{Constructing the risk estimators for mixture model selection}
Till now we have focused on the form and the smoothing parameter of the kernel for our problem.
Using this kernel we will now illustrate the steps in estimating  the   risk of each parametric mixture model $\mathcal M_k,~~k=1,2,\ldots$ along
with the risk of the empirical model. First, we compute the estimates of the
parameters $\theta_k$ for each parametric model $\mathcal M_k$. Simultaneously, we construct
a projection matrix $\mathbb{P^\ast_{\hat\theta}}$ for each model, which is strictly based on the
estimates ${\hat{\theta}_k}$ and the observations. We also
calculate the score-centered kernel  $\mathbb{K}_{scen({M_{\hat{\theta}_k}})}$. Then we calculate the
risk estimate of each model using \eqref{eq:riskfinal} and \eqref{eq:qbicrisk}.

\section{ Simulation results and Applications}
\label{sec:simul}

Now we apply our model selection tool to select the number of components in
mixtures. The following steps were used
for choosing models in all the examples of this section. First,we
standardized the data (i.e.  the variance of each
variable was scaled to unity). This step was done because we used the same $h$ for all the
variables. In this simulation study we restrict ourselves to normal kernels based on the discussion in \ref{sec-kernel-specification}.
 Using the results in Section \ref{sec:sdof}  we estimated the $sDOF$,
which  gives us a range of interesting smoothing parameter
values. The simulation results reported in this section 
 used one representative value from this range, though the results were 
 quite stable across our suggested range.
 The ML estimates for the range of models under consideration
are calculated using an EM Algorithm and the projection matrices $\mathbb
P^\ast_{\hat \theta}$  based on this estimate. Note that for analysis at
different levels of smoothing  we only need to recalculate the 
  distance and
computationally simple  matrix
$\mathbb K$,  and not the projection matrix   $\mathbb P^\ast_{\hat \theta}$.

\subsection{ Simulation experiment} 
\begin{figure}[tbp]
  \centering 
\includegraphics*[angle=270,width=.24\textwidth, trim= 0 0 0 0]{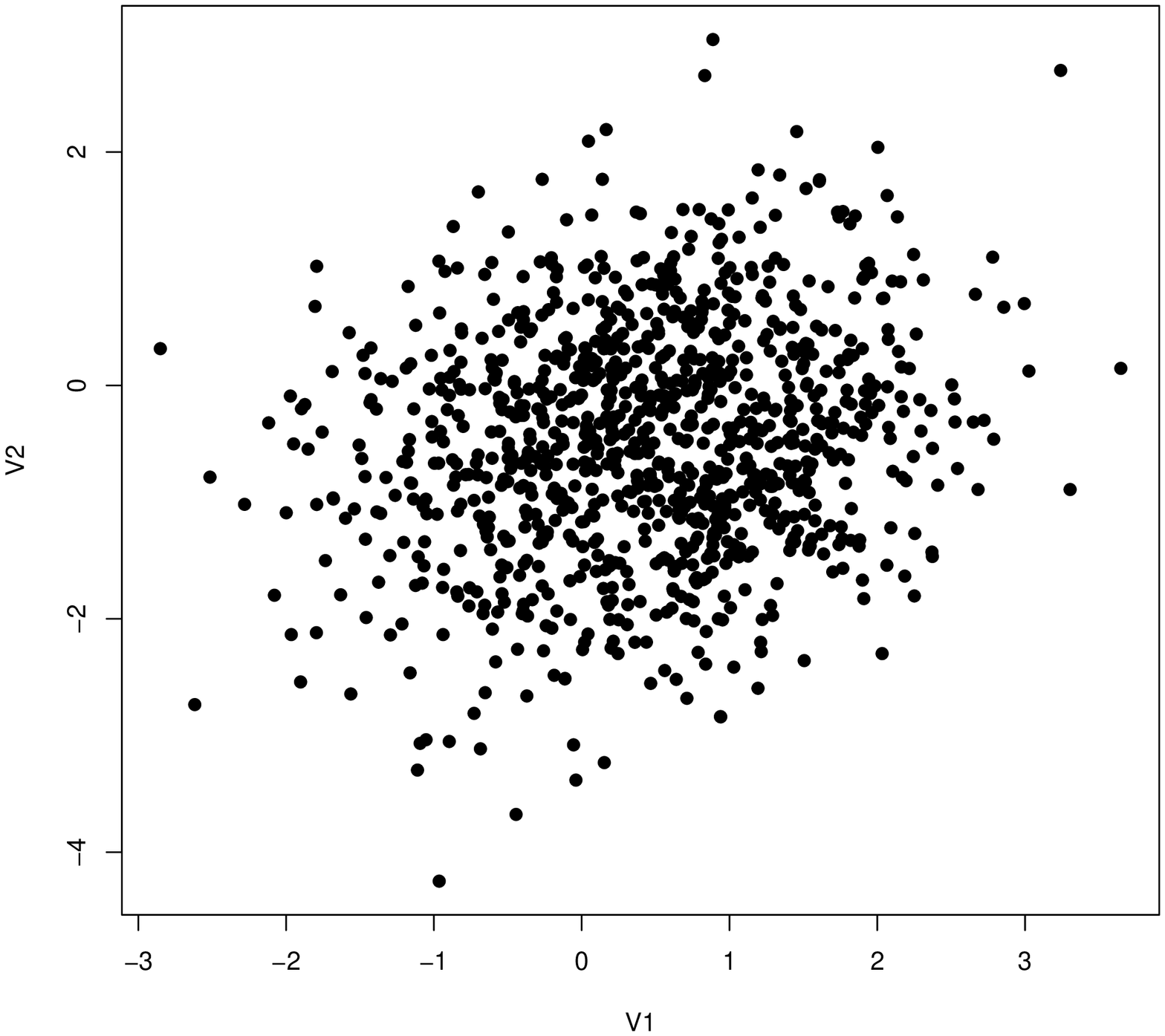}
\includegraphics*[angle=270,width=.24\textwidth, trim= 0 0 0 0]{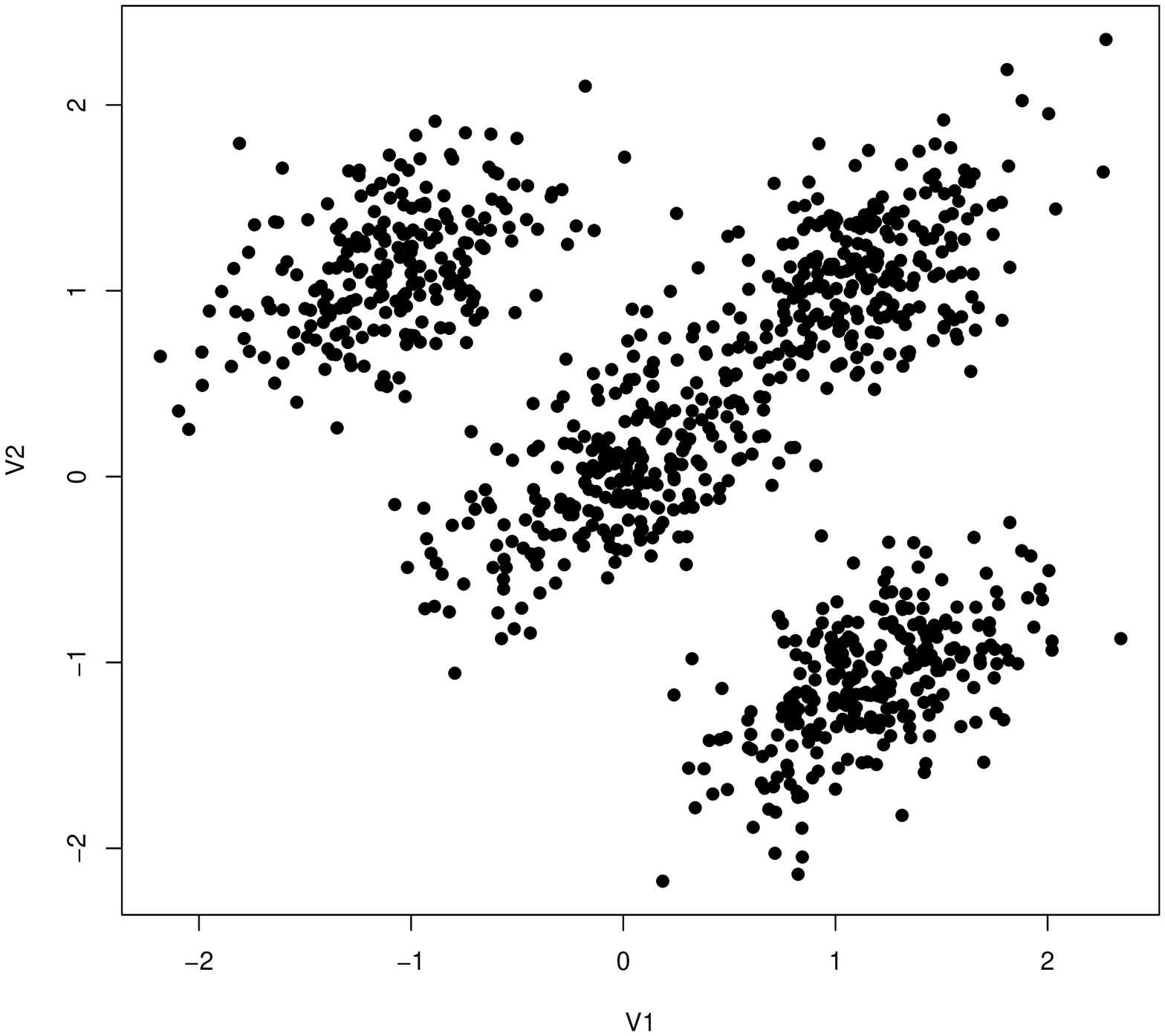}
\includegraphics*[angle=270,width=.24\textwidth, trim= 0 0 0 0]{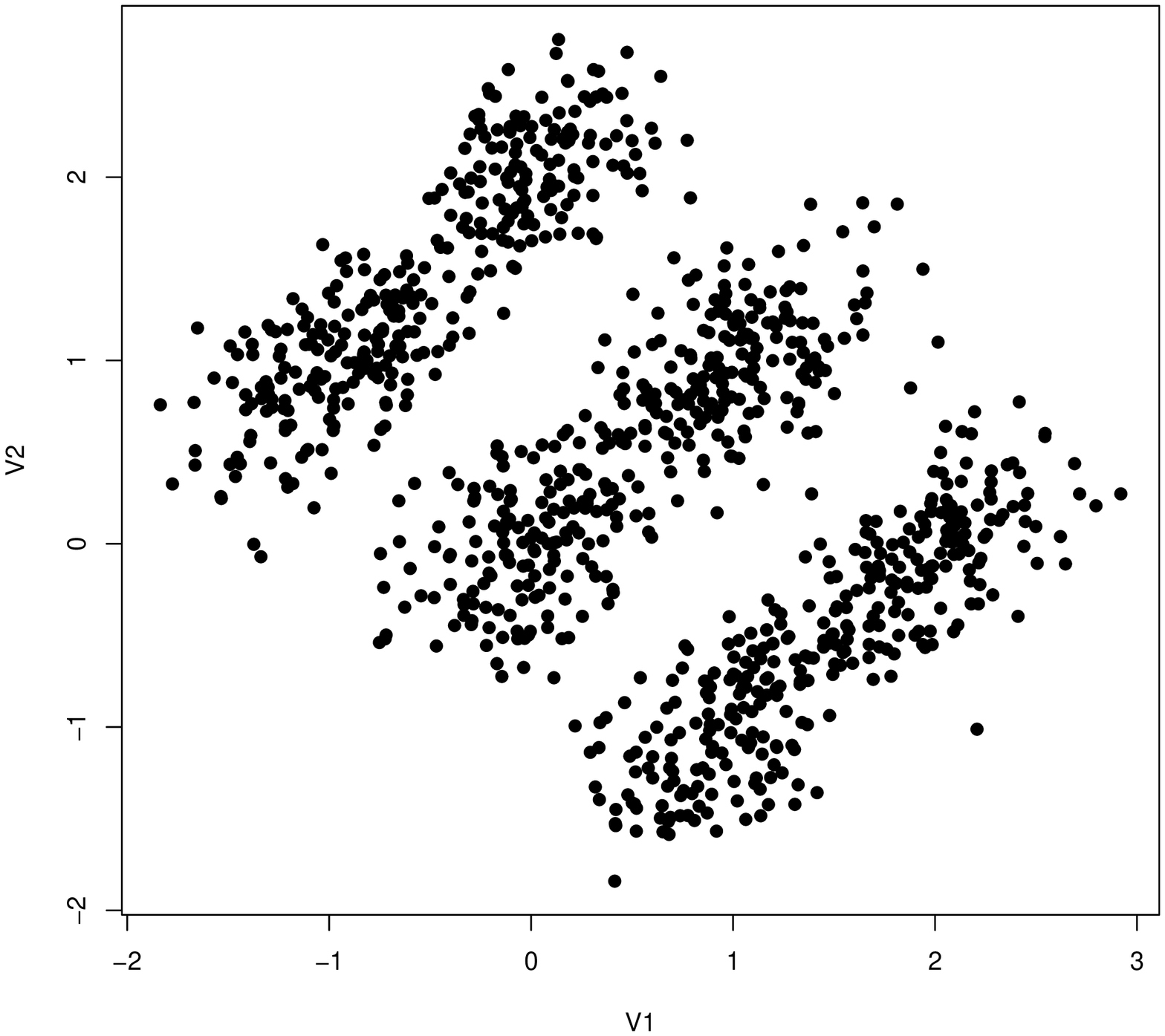}
\includegraphics*[angle=270,width=.24\textwidth, trim= 0 0 0 0]{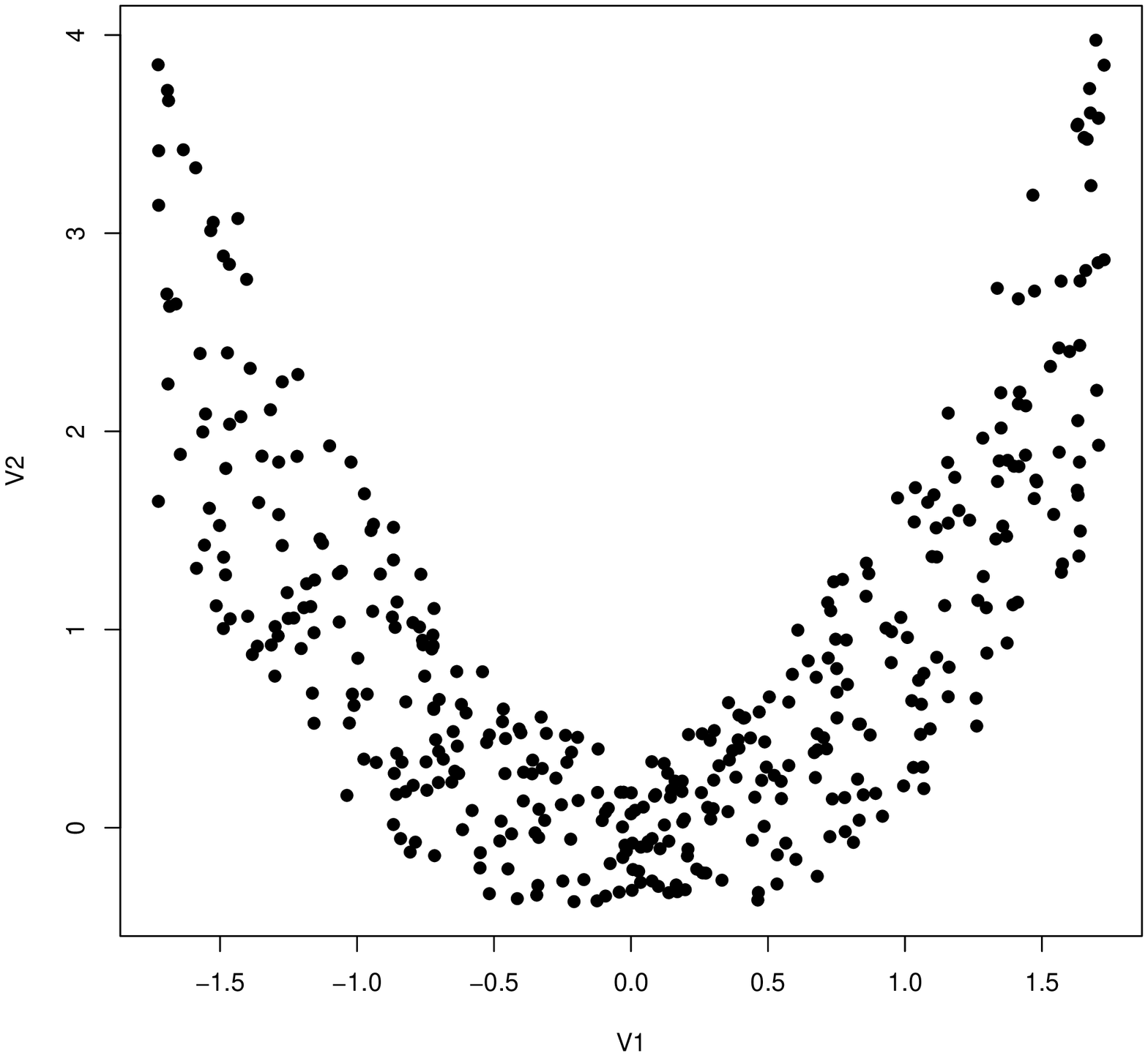}

\noindent
~\hfill Model (1) \hfill \hfill Model (2) \hfill \hfill Model (3) \hfill \hfill Model (4) \hfill~
\caption{Representative datasets from two dimensional simulation}
\label{fig:dataset}
\end{figure}

Given a mixture parameter set and a fixed sample size,
 we generated 100 sets of random samples from the true mixture
 distribution. 
We have tabulated
 the frequency of selection of the number of
 components  by each of the 5
 model selection methods (AIC, QAIC, MRA, BIC,  and QBIC) under  different
 examples, in an experimental design that varies by dimension,
 sample size and separateness of the components. For some example we also vary the range of models in the model space.

Here, we have fixed the variance structure, so our model selection only refers to
selecting  the number of components. Readers familiar with the MCLUST
software
 should note that \citet{mclust} provide a unified
strategy to select both the number of components and the variance structure,
but only with a single criterion, the BIC.  We now discuss the results of our simulation experiment.




\begin{table}
\caption{\label{tbl:2dim}Model Selection results for a set of 2 dimensional
  simulated datasets. \t{Boldface} indicates selection of true model.}
\centering
\begin{small}
\begin{tabular}{>{\columncolor[gray]{0.95}}r r r r r r r r r r r  }
\hline
\hline
 \rowcolor{white}~ &\multicolumn{10}{c}{ Estimated number of Components}\\
\cline{2-11}
\rowcolor[gray]{0.95} & 1  & 2 & 3 & 4 & 5 &  6 &  7 &  8 &  9 & $>$10 \\
\simtitle{Model (1): 2 moderately overlapping clusters,2D,  n=1000}
AIC   & \0 &\t{ 63} & 16 & 6 & 4 & 4 & 3 & 2 & 1 & 1 \\
\RAIC & 5 & \t{ 62} & 19 & 2 & 3 & 6 & 3 & \0 & \0 & \0 \\
\RADE   & 100&\t{ \0} & \0 & \0 & \0 & \0 & \0 & \0 & \0 & \0 \\
BIC   & 48 &\t{ 52} & \0 & \0 & \0 & \0 & \0 & \0 & \0 & \0 \\
\RBIC & 21 &\t{ 77} & 2 & \0 & \0 & \0 & \0 & \0 & \0 & \0 \\
\simtitle{Model (2): 4 distinct clusters, 2D, n=1000}
AIC   & \0 & \0 & \0 & \t{52 } & 17 & 8 & 12 & 3 & 3 & 5 \\
\RAIC & \0 & \0 & \0 & \t{66 }& 18 & 6 & 4 & 4 & 1 & 1 \\
\RADE   & \0 & \0 & \0 & \t{100} & \0 & \0 & \0 & \0 & \0 & \0 \\
BIC   & \0 & \0 & \0 & \t{100} & \0 & \0 & \0 & \0 & \0 & \0 \\
\RBIC & \0 & \0 & \0 & \t{98 }& 2 & \0 & \0 & \0 & \0 & \0 \\
%
%
\simtitle{Model (3): 6 distinct clusters, 2D, n=1000}
%
AIC   & \0 & \0 & \0 & \0 & \0 & \t{43 }& 34 & 12 & 4 & 7 \\
\RAIC & \0 & \0 & \0 & \0 & \0 & \t{64 }& 31 & 9 & 4 & 2 \\
\RADE   & \0 & \0 & \0 & \0 & \0 & \t{100} & \0 & \0 & \0 & \0 \\
BIC   & \0 & \0 & 5 & \0 & \0 & \t{95} & \0 & \0 & \0 & \0 \\
\RBIC & \0 & \0 & 1 & \0 & \0 & \t{94 }& 5 & \0 & \0 & \0 \\
%
\simtitle{Model (3): 6 distinct clusters, 2D, n=300}
%
%
AIC   & \0 & \0 & \0 & \0 & \0& \t{1 }& 1 & 6 & 28 & 64 \\
\RAIC & \0 & \0 & \0 & \0 & 1 & \t{28} & 38 & 17 & 10 & 6 \\
\RADE   & \0 & 3 & 25 & 16 & 13 & \t{38} & 5 & \0 & \0 & \0 \\
BIC   & \0 & \0 & 50 & 26 & 12& \t{12 }& \0 & \0 & \0 & \0 \\
\RBIC & \0 & \0 & 11 & 5 & 17 &  \t{49} & 17 & 1  & \0 & \0 \\
\hline
\end{tabular}
\end{small}
\end{table}

\subsubsection{Two dimensional simulation}
In our first study we held the data dimension to two and varied the normal
mixture model as well as the sample size (for some cases). The scatterplot of a representative sample from the simulation of each of the 4 models is given in Figure 2. 

First we consider  simulations from the first 3 models which are generated from 2, 4 and 6 component normal mixtures respectively.
From Table~\ref{tbl:2dim}  it is clear that AIC and QAIC are the more liberal methods-
in the sense of favoring models larger than the true model.  Although QBIC and BIC are both
potentially conservative, QBIC is slightly less so and made model selection
errors on both sides of the true value. The MRA criterion (here using the $m=n$ risk) was fairly
conservative in some cases (e.g. Simulation 1) but otherwise fairly competitive. 
 Like BIC it never overshot the true model.

Next we explore the effect of sample size on model selection in Model (3).
Decreasing the sample size from 1000 to 300 tended to increase the
underestimation of the three conservative methods, but did not greatly degrade
the liberal estimates of AIC and QAIC, reflecting their lack of consistency.

Overall QAIC performed much better than AIC. Though QAIC overestimated the
  complexity, in comparison to AIC the distribution of the the number
  of components was much more peaked at the true number while the right tail
  tapered off very quickly. Also, if the components had high degree of
  overlap,  in the sense of not displaying distinct modes for distinct
  components, then QBIC performed better than BIC, which most
  often underestimated the number of components.

Note also that the quadratic risk based methods are more stable with change
of sample size. In those cases where the sample size was small and BIC badly
underestimated, the QBIC and MRA methods were distinctly better.
Finally, the simulation study revealed that
 adequacy is a very useful measure. In many examples it provided
 extra information about the minimal number of components that  would be needed to explain
  the data in hand. Often the minimal adequate model had a considerably smaller  number of
  components than the one with minimum risk QAIC model.

\subsection{Quadratic risk as global measure of risk}
Next we explore two examples where the true model is beyond our model space and show how the  quadratic risk measures provides an excellent  global model selection criterion in identifying the scenario.
  We have already observed that with n=1000 in Model (3) with six distinct normal components, QAIC, QBIC, MRA and BIC all overwhelmingly select  six components. But what happens if we explore only up to 5 components, i.e we do not include the true model in the model class?  Table 2 summarizes the results for the 5 different criteria when we explore only up to 5 components. 
We immediately  observe that BIC selects the 3 component model in 100\% of
  the cases. Independently our quadratic measures  also selected this model in 52\%  and 87\% of the cases, but the MRA criterion shows that none of the models are adequate.
 Thus if we based our conclusion on BIC alone we would have failed to recognize that we are trapped in a local minima. But if we use the quadratic risk measure along with the MRA criterion we are forced to explore models beyond 3 components as we have a clear indication that even the 5 component model is not adequate.

\begin{table}
\caption{\label{tbl:global} Model selection results demonstrating global risk measure}
\centering
\begin{small}
\begin{tabular}{>{\columncolor[gray]{0.95}}r r r r r r rrrrr }
\hline
\simtitle{Model (3): 6 components, 2D, n=1000 (limited to 5 components)
 }
 \rowcolor{white}~ &\multicolumn{10}{c}{ Estimated number of Components}\\
\cline{2-11}
& 1  & 2 & 3 & 4 & $>$ 5\\ \hline
  AIC & \0   & \0   & 52   & 9   & 39   \\
\RAIC & \0   & \0   & 58   & 6   & 36   \\
\RADE & \0   & \0   & \0   & \0   & \0   \\
  BIC & \0   & \0   & 100   & \0   & \0   \\
\RBIC & \0   & \0   & 87   & \0   & 13   \\
\hline
\end{tabular}
\end{small}
\end{table} 

We now explore an example where the model is not a normal mixture, but
rather  when the scatter of the distribution has
the shape of the letter U, as given in Figure 2(Model 4). 
These datasets were simulated by first generating  the $x$ coordinate uniformly between
-1.5 and 1.5 and then generating $y$ uniformly between $x^2-1$ and $x^2+1$.
By design, this density is very hard to capture with a few normal components. So we
 explored up to 14 normal components. BIC and QBIC mostly choose 8-10
 components. Even using the liberal measure QAIC, 51\% of the times we choose
 a model with less than 14 components. But the MRA criterion overwhelmingly (70\% of the cases) shows
that even a 14 component model is not adequate. This suggests that we should either explore beyond 14 component model or we should explore beyond normal mixtures.

\begin{table}
\caption{\label{tbl:U} Model selection results for the U Example (simulation (4))}
\centering
\begin{small}
\begin{tabular}{>{\columncolor[gray]{0.95}}r r r r r r r r r r r r r r r}
\hline
 \rowcolor{white}~ &\multicolumn{10}{c}{ Estimated number of Components}\\
\cline{2-11}
\rowcolor[gray]{0.95}  & 1 & 2 & 3 & 4 & 5 & 6 & 7 & 8 & 9 & 10 & 11 & 12 & 13 & $>$14 \\
  \hline
AIC &\0&\0&\0&\0&\0&\0&\0&\0&\0&\0&\0& 2 & 16 & 82 \\
\RAIC &\0&\0&\0&\0&\0&\0&\0&\0& 1 & 5 & 7 & 19 & 19 & 49 \\
 \RADE &\0&\0&\0&\0&\0&\0&\0&\0& 1 & 3 & 7 & 8 & 3 & 8 \\
  BIC &\0&\0&\0&\0&\0&\0&\0& 16 & 50 & 24 & 6 & 3 &\0& 1 \\
  \RBIC &\0&\0&\0&\0&\0&\0&\0& 12 & 30 & 28 & 9 & 10 & 7 & 4 \\
   \hline
\end{tabular}
\end{small}
\end{table}

In both these examples it can be easily seen  that BIC, which is the most widely used model selection method for mixtures, does not provide the best result, simply because it is  a relative rather than a global measure of risk. On the other hand in each of these simulations the MRA criterion clearly points out, that the minimum risk model does not lie in the model class that we are considering.

\begin{table}
\caption{High dimensional simulation results}
\begin{small}
\begin{tabular}{>{\columncolor[gray]{0.95}}r r r r r r r r r r r || r r r r r p{.5in}}
\hline
 \rowcolor{white}~ &\multicolumn{10}{c||}{ Estimated number of Components} & \multicolumn{6}{c}{ Estimated number of Components}\\
 \rowcolor{white}~ &\multicolumn{10}{c||}{ } & \multicolumn{6}{c}{ (search limited to 5 components)}\\
\cline{2-17}
\rowcolor[gray]{0.95} & 1  & 2 & 3 & 4 & 5 &  6 &  7 &  8 &  9 & $>$10  & 1  & 2 & 3 & 4 & $>$5 &~ \\
\simtitleg{Model (5):  6 components, n=1000 in 4 dimensions}
AIC & \0 & \0 & \0 & \0 & \0 & \t{23} & 40 & 14 & 15 & 8 &       \0 & \0 & \0 & \0 & 100 \\
  \RAIC & \0 & \0 & \0 & \0 & \0 & \t{56} & 27 & 11 & 3 & 3 &    \0 & \0 & 8 & \0 & 92 \\ 
  \RADE & \0 & \0 & \0 & \0 & 4 & \t{96} & \0 & \0 & \0 & \0 & \0 & \0 & \0 & \0 & 4 \\
  BIC & \0 & \0 & 20 & \0 & 3 & \t{77} & \0 & \0 & \0 & \0 &   \0 & \0 & 88 & \0 & 12 \\
  \RBIC & \0 & \0 & 3 & \0 & 3 & \t{95} & \0 & \0 & \0 & \0 & \0 & \0 & 49 & \0 & 51 \\
\simtitleg{Model (6):  6 components, n=1000 in 8 dimensions}
AIC & \0 & \0 & \0 & \0 & \0 & \t{17} & 32 & 22 & 16 & 13 &     \0 & \0 & \0 & \0 & 100 \\
  \RAIC & \0 & \0 & \0 & \0 & \0 & \t{54} & 34 & 10 & 2 & \0 &  \0 & \0 & 33 & 1 & 66 \\ 
  \RADE & \0 & \0 & \0 & \0 & 8 & \t{92} & \0 & \0 & \0 & \0 &\0 & \0 & \0 & \0 & 8 \\
  BIC & \0 & \0 & 25 & \0 & 5 & \t{69} & 1 & \0 & \0 & \0 &  \0 & \0 & 71 & 0 & 29 \\
  \RBIC & \0 & \0 & 5 & \0 & 4 & \t{88} & 2 & \0 & \0 & \0 &\0 & \0 & 84 & \0 & 16 \\
\simtitleg{Model (7):  6 components, n=1000 in 12 dimensions}
    AIC & \0 & \0 & \0 & \0 & \0 & \t{12} & 31 & 18 & 16 & 23 &  \0 & \0 & 13 & 7 & 80 \\ 
  \RAIC & \0 & \0 & \0 & \0 & \0 & \t{49} & 34 & 9 & 5 & 3 &    \0 & \0 & 75 & 21 & 4 \\ 
  \RADE & \0 & \0 & \0 & \0 & 10 & \t{88} & 1 & 1 & \0 & \0 &  \0 & \0 & \0 & \0 & 10 \\
    BIC & \0 & \0 & 27 & \0 & 4 &  \t{63} & 5 & 1 & \0 & \0 &  \0 & \0 & 93 & 3 & 4 \\  
  \RBIC & \0 & \0 & 7 & \0 & 3 &  \t{81} & 6 & 2 & \0 & \0 &  \0 & \0 & 90 & 7 & 3 \\  
\hline
\end{tabular}
\end{small}
\end{table}

\subsubsection{High dimensional simulation}
We now examine the effect of increasing dimensions.  Our 4D, 8D and 12D
  examples  with 6 components  were generated  in the following way: The
  first two dimensions are  the same as 2D Model (3) used in Section
  5.1.1; the remaining dimensions are generated from standard normal
  variates.
Like the 2D Model (3) the performance of the 5 criteria  are judged under  two  scenarios (i) when the true model is in the model space   (ii) and when the true model was not included in the search space. The results for both scenarios  are given in Table 4, the left panel corresponding to scenario (i), where we explore up to 10  components and the right panel corresponding to scenario(2) where we restrict ourselves to less than 6 components. 

 First we summarize the results from scenario (i). AIC grossly overestimates the number of components.
 QAIC overestimated the number of components but less so than the AIC, but
 more importantly its performance did not deteriorate  much with the increase
 of dimensions. As  dimension increases, BIC tends to  underestimate the true
 number of components with higher percentage, and in many cases (20\%, 23\%,
 25\% for 4, 8 and 12 dimensions respectively) BIC even chose the 3 component
 model. QBIC and MRA erred on either side of the true model, but both of them
 chose the true model with higher percentage than BIC. Moreover, when QBIC
 and MRA  underestimated they erred slightly selecting 5 components, while BIC selected 3 component model in many cases.

For scenario (ii) the results  are consistent with the findings of the two
dimensional case, i.e. the various information criterion overwhelmingly suggest that the 3 component model is the true model, whereas the MRA criterion combined with the quadratic measures clearly points out (in  96\%, 92\%, 90\% for 4, 8 and 12 dimensions respectively) that we have not reached a global minimum and we should expand our search space.

These examples clearly show the usefulness of quadratic risk based methods in high dimensions, both for selecting the true model and providing a global infrastructure for model selection.

\subsection{Application to real dataset}

We  apply our model selection criterion to a  real dataset where the goal is to identify the number of gene groups displaying distinctive gene expression profile under a set of conditions.
 The details of the experiment is described in \citet{Chit:2002} and our
 dataset was  obtained through personal communication from Dr. F. Pugh of Pennsylvania State University.  Here we present a very short description of the experimental design and preprocessing of the data.

The experiment was designed to answer an interesting  questions in transcriptional genomics: Is the binding phenomenon of TATA binding proteins (TBP)  affected by neighboring TAF1 proteins? The TBPs bind to the TATA box and play a crucial role in transcription which finally results in protein synthesis. The TBP is intertwined to form a 3D structure in such a way that it forms a concave surface and a convex surface;
the concave surface attaches to the TATA box, whereas the convex surface
may attach to some other protein, for example the TAND region of  TAF1 proteins. To understand this interplay of TBP and the TAND region of TAF1, a set of 19 experiments, were designed using 6226 yeast genes for each strain. The two main factors that were varied in these experiment are (i) whether the TAF1 proteins of the yeast cells had the the TAND regions intact (wild types) or whether the TAND regions were removed ($\Delta$-TAND's), (ii) whether the concave region of the TBP proteins remained unmutated (wild type) or whether a yeast strain with mutation was used. Further 3 different types of mutations  and a few controls were also used. Gene clusters displaying distinct gene expression profiles, one  for each of the 4 combinations of the two levels of the two major factors would suggest a ``interaction'' model, whereas only two clusters would suggest that TBP's are not influenced by neighboring proteins.

The analysis was done as follows. After leaving out genes which did not show
significant change over different conditions, we analyzed n=2809 genes under
d=19 conditions (dimensions), completely ignoring the information on the identity of the conditions. We used normal mixtures to fit the data and explored up to 8 components.
The risk curves of each of the 4 model selection criteria is given in Figure~\ref{fig:risk}. Based on the minimum risk,  BIC chose three components where as AIC chose nine components.
For calculating the quadratic risk we first applied the sDOF analysis detailed in Section 4.2.2 and chose a range of $h$. The risk curves given in Figure~\ref{fig:risk} are for $h=1.7$ for which the sDOF was  226, but the minimum risk decision remained stable for a wide range of smoothing parameters. Based on the quadratic risk calculations QAIC chose six components and QBIC chose four components. Moreover using the MRA criterion it was clear that four components was adequate.

There is a sense in which selecting four components is the right answer for this problem. 
Biologists conjecture that if the interaction model is true there should be 4
groups of genes each displaying a separate profile under the 4
conditions. Among the 3 clusters given by BIC, one 
 profile matched with the distinct expression profile for $\Delta$-TAND with
no  mutation, but the other two gene clusters were mixed.
On the other hand the four clusters given by QBIC clearly identifies four distinct
patterns under the major combinations, thereby providing a meaningful summary
and cluster analysis of this high dimensional data. Additionally QAIC and
more prominently QBIC shows a local minima for a two components mixture. This
suggested there might be two big groups. Further exploration shows that the
two top clusters differentiated the mutated and unmutated TBP yeast strains.

\begin{figure}[tbp]
  \centering 
\includegraphics*[width=.7\textwidth,height=3in, trim= .2in .2in .1in .2in]{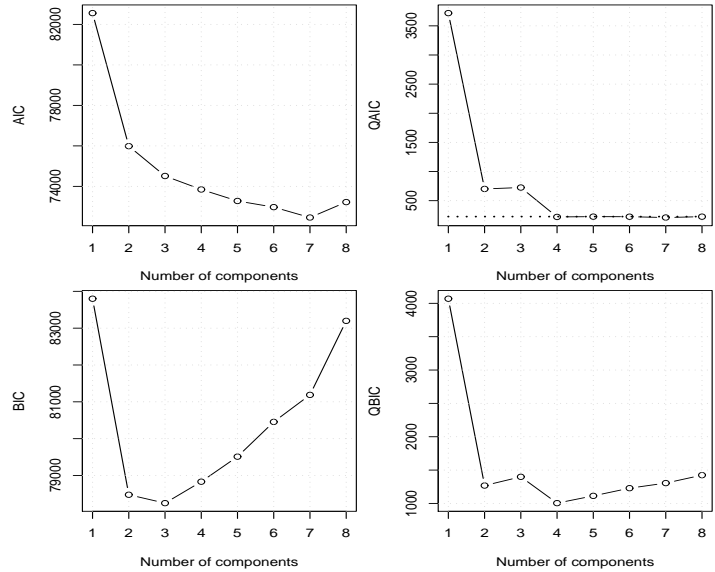}
\caption{Risk Calculation for the yeast gene expression dataset}
\label{fig:risk}
\end{figure}

\section{Conclusion and Future Direction} \label{sec:discussion}

In this paper we  have developed a comprehensive tool for high dimensional model selection,
using quadratic risk. Key to the calculation of quadratic risk is its representation
and decomposition using appropriate projection operators, and the 
estimation of these operators using empirical versions.  Decomposing  the quadratic risk into the model
lack of fit and the parameter estimation cost 
enabled us to build a method, QBIC, that mimics BIC.

One feature of the derivation of our methodology, one that separates it from the AIC and BIC, is that the asymptotic expansions are based on nonparametric goodness of fit tests, not likelihood ratio asymptotics.  We believe that our unified approach for building and estimating quadratic risk
could pave the way for designing appropriate model selection criterion for a
host of previously unsolved problems, especially where the irregularity of
parameter space eliminates the asymptotic theory underlying the use of AIC, BIC and similar other
criteria. In addition, the risk of the empirical estimate, which is a
natural outcome of our risk analysis, provides a threshold that enables us to
assess whether the model that was selected was itself a good fit.  
Unlike AIC and  BIC where only the
optimum model is chosen, model adequacy provides us with the extra
information of whether any of the models are adequate, or if there is a range of
models that are adequate.

In addition, we  also showed that the use of an appropriate kernel could
enable one to   minimize the computational burden  associated with a
high-dimensional problem. Moreover,  there is certainly  more to learn about the structure of
a dataset than  can be revealed by analyzing a data with a single model selection
criterion.  This is an area of future research.

\small
\bibliography{mixture,model}
\bibliographystyle{Chicago}

\section*{ Appendix: Proof of Proposition \ref{pro:pec}} \label{app:proof}}
 \begin{appendix}

  We start with an asymptotic approximation to $\mbox{Parameter Estimation
Error}=\left[ d_{K}(F_{\tau },M_{\hat{\theta}})-d_{K}(F_{\tau },M_{{\theta }%
_{\tau }})\right] $. Now we define a function in the $\theta $-space, by $%
L(\theta _{\tau },\hat{\theta})=d_{K}(M_{\hat{\theta}},M_{\theta _{\tau }})$%
, which has the form 
\begin{equation*}
L(\theta _{\tau },\hat{\theta})=\iint K\left( x,y\right) \ d(M_{\hat{\theta}%
}-M_{\theta _{\tau }})(x)d(M_{\hat{\theta}}-M_{\theta _{\tau }})(y).
\end{equation*}%
For  density function $m_{\hat{\theta}}$ we assume that we can create a
Taylor expansion in $\theta _{\tau }\ $ such that 
\begin{equation*}
(m_{\hat{\theta}}-m_{\theta _{\tau }})(x)=(\hat{\theta}-\theta _{\tau })^{T}%
\mathbf{u}_{\theta _{\tau }}(x)\ dM_{\theta _{\tau }}(x)+o\left( |\hat{\theta%
}-\theta _{\tau }|\right) .
\end{equation*}
Thus we get the following quadratic approximation for the loss 
\begin{equation}
mL(\theta _{\tau },\hat{\theta})\!=\!m(\hat{\theta}-\theta _{\tau })^{T}\left[
\iint \!\!\mathbf{u}_{\theta _{\tau }}(x)K(x,y)\mathbf{u}_{\theta _{\tau
}}(y)^{T}\ M_{\theta _{\tau }}(x)M_{\theta }(y)\right] (\hat{\theta}-\theta
_{\tau })\!+\!o(m|\hat{\theta}-\theta _{\tau }|^{2})  \label{eq:quadloss}
\end{equation}%
Next, under standard regularity assumptions, it is a well known result for
maximum likelihood that 
\begin{equation*}
m^{1/2}(\hat{\theta}-{\theta _{\tau }})=m^{1/2}\mathbb{J}^{-1}\left( \frac{1%
}{m}\sum_{i}\mathbf{u}_{\theta _{\tau }}(x_{i})\right) +o_{p}(1).
\end{equation*}%
This shows that the error term in \eqref{eq:quadloss} is $o_{p}(1).$ In
addition, substituting this expression into the first term gives the
relationship: 
\begin{equation}
mL(\theta _{\tau },\hat{\theta})=m\iint V(x,y)\ \hat{F}(x)\hat{F}(y)\
+o_{p}(1)
\end{equation}%
as needed for the proposition. Note that $\int V(x,y)M_{\theta _{\tau }}(y)=0
$, so the asymptotic mean for $mL(\theta _{\tau },\hat{\theta})$  is  
\begin{equation*}
E\left[ \iint V(x,y)\ \hat{F}(x)\hat{F}(y)\right] =trace_{M_{\theta _{\tau
}}}(V).
\end{equation*}

 \end{appendix}
 \end{document}